\documentclass[12pt]{article}
\usepackage[koi8-u,cp1251]{inputenc}
\usepackage[english,russian,ukrainian]{babel}
\usepackage{citehack,amsmath,amsthm,amssymb}
\usepackage{psfrag}
\usepackage[active]{srcltx}
\usepackage{graphicx}
\usepackage{epsfig}
\usepackage{enumerate}
\usepackage[a4paper,left=28mm,right=12mm,top=2cm,bottom=2cm]{geometry}

\usepackage{xspace,cite,indentfirst}
\usepackage[mathscr]{eucal}
\theoremstyle{plain}
\theoremstyle{definition}

\begin{document}

\begin{center}
{\bf ROUGHNESS OF DICHOTOMY FOR THE INTERCONNECTED SYSTEM OF OPERATOR-DIFFERENTIAL EQUATIONS}
\end{center}

\begin{center}
{\bf Iskra O. (Institute of mathematics of NAS of Ukraine), Dashkovskiy S. (University of W$\ddot{u}$rzburg), Pokutniy A. (Institute of mathematics of NAS of Ukraine), Slynko V. (University of W$\ddot{u}$rzburg, ),  Zahorulko Yu. (Institute of mathematics of NAS of Ukraine)}
\end{center}

{\bf Introduction.} The concept of exponential dichotomy and its connection with nonlinear dynamics takes inspiration from Lyapunov's paper \cite{Lyap} (1892) on stable manifolds. This fact was known and used by Poincare \cite{Poin}, \cite{Poin1} (1890, 1892) in the study of transversal intersections for stable and unstable manifolds. It should be noted that problems of the existence of bounded solutions appeared in the study of the question of the existence of periodic solutions for the linear and nonlinear systems in the works of the above-mentioned scientists. The condition of the existence of bounded solutions of a nonhomogeneous equation
\begin{equation} \label{eq:1}
x'(t) = A(t)x(t) + f(t),
\end{equation}
 with bounded functions $f(t)$ for finite-dimensional systems of linear equations with variable coefficients and continuous time appears in the works of Perron \cite{Perr} -- \cite{Perr1} (1928, 1930) in the study of nonlinear perturbations of such equations (that generalizes  Hadamard's results \cite{Hadamar} (1901)) and in the paper of his student Li Ta \cite{Ta} (1934) for systems with discrete time (in these papers the notion of dichotomy was not yet defined). A similar problem for a nonlinear equation with a stationary linear part was considered by Bohl. For bounded coefficient $A(t)$ the dichotomy condition and it's the equivalence to the existence of bounded solution condition of the corresponding nonhomogeneous equation appeared and established in the paper of Meisel \cite{Mais} (1954), who used the technique of the Perron transformation. For the nonlinear systems of differential equations, this theory began applied after the well-known Mel'nikov paper \cite{Meln} (1964).
The generalization of these results to the case of infinite-dimensional spaces with bounded
operator coefficients was carried out by many mathematicians, in particular, in the papers of
famous mathematicians Massera, Schaeffer \cite{Massera1}, \cite{Massera2} (1958, 1966), Daletskii, Krein \cite{dal} (1970), Sacker, Sell  (1976, 1978, 1979), Hirsch, Pugh and Shub \cite{Hirsch} (1977), Pliss
\cite{Pliss} - \cite{Pliss3} (1964, 1966, 1969, 1977), Shen and Yee \cite{Shen} (1966). Note that in the papers of Massera, Schaeffer one of the equivalent definitions of dichotomy (for certain classes of equations) was considered, namely - the notion of admissibility of functional spaces $(B, D)$ in which the equation (\ref{eq:1}) was investigated, where the nonhomogeneity $f(t) \in D$ and solution $x(t) \in B$ respectively. A pair $(B, D)$ is admissible if for each $f(t) \in D$ there is a solution $x(t)$ of the nonhomogeneous equation (\ref{eq:1}) such that $x(t) \in D$. Concept of admissibility is like similar to the
concept of stability under constantly acting perturbations. The first to draw attention to the 
importance of these properties was Perron \cite{Per}, who showed that the admissibility of the pair space $(C, C)$ ($C$ is the space of bounded and continuous functions on the semi-axe $\mathbb{R}_{+}$ with values in some space $X$) under sufficiently general assumptions is equivalent to the existence of an exponential dichotomy. Hartman \cite{Hartman} (1964) contributed to the development of such questions in the study of higher-order differential equations. A series of results in this direction is summarized and generalized in Coppel's papers \cite{Coppel1}, \cite{Coppel} (1978).
One of the applications of the concept of exponential dichotomy in finite-dimensional and infinite-dimensional spaces is the connection with the Floquet theory for time-dependent (or almost periodic) linear systems, which was developed in the works of Johnson \cite{Johnson}, \cite{Johnson1},
Moser, Johnson \cite{Johnson}, \cite{Johnson2}, Chow, Malle-Paret \cite{Chow}, \cite{Chow1}. Of particular interest are the investigations of Anosov and Sinai \cite{Anosov}, \cite{Anosov1} where the corresponding definitions were used in the study of flows and cascades on a smooth manifold.

Another impetus in the study of the problem of solutions bounded on the whole axis was given by the Palmer in the his well-known lemma, which connects the concept of exponential dichotomy on the semi-axes $\mathbb{R}_{+} = [0; +\infty)$, $\mathbb{R}_{-} = (-\infty, 0]$ with Noethirianity (Fredholm with nonzero index) of the corresponding differential operator. In 1984, he proved \cite{Palmers} the assertion that if
a linear system
\begin{equation} \label{eq:2}
x'(t) = A(t)x(t)
\end{equation}
admits an exponential dichotomy on the semiaxes $\mathbb{R}_{+}$, $\mathbb{R}_{-}$
then the operator which acts according
to the rule
\begin{equation} \label{eq:3}
(Lx)(t): = x'(t) - A(t)x(t)
\end{equation}
is Noetherian. Later \cite{Palm} 1988 Palmer proved the opposite statement for bounded and continuous on the interval $J$ matrix-valued function $A(t)$. He proved that if $L: BC^{1}(J) \rightarrow BC(J)$ is semi-Fredholm, then the homogeneous system (\ref{eq:2}) admits an exponential dichotomy on $J$.
Related results were obtained by Ben-Artzi and Hogberg \cite{Ben} (1992) on the space of integrable functions. Earlier, this result was proved in one direction by Sacker \cite{Sackk} (1978). Note that this idea was further developed in the works of Samoilenko and Boichuk \cite{BoiSam1} (2016), where with using of generalized inverse and pseudoinverse Moore-Penrose matrices, the considered porblem of the existence of bounded on the whole axis solutions under linear and nonlinear perturbations was investigated. For equations in Banach spaces with bounded coefficient $A(t)$ this result was obtain by  Boichuk, Pokutnyi \cite{Boi_Pok} (2006), and in the paper \cite{Pok1} (2012) for equations
with an unbounded coefficient under nonlinear perturbations. The Fredholm property of the
corresponding differential equation with unbounded operator coefficients is studied in the
paper of Latushkin, Tomilov \cite{Lat} (2005) for the abstract parabolic operators. Similar statements for the dichotomy on the semi-axes $\mathbb{R}_{+}, \mathbb{R}_{-}$ was established by Rodriguez and Ruas-Filho \cite{ROdr} (1995). They proved an analogue of Fredholm alternative. The same assertions for difference equations in Banach spaces were studied in the works of Chueshov \cite{Chuesh} (2002), Slyusarchuk \cite{Slyus} (2003), Bareira, Valls \cite{Bar} (2011), Bento, Silva \cite{Bento} (2012), Samoilenko, Boichuk \cite{BoiSam1} (2016). The generalization of these results to more general cases in finite-dimensional and infinite-dimensional spaces remains actual today. Our goal is to further study the theory of exponential dichotomies for differential equations with unbounded operators in the linear part.

The concept of exponential dichotomy for evolutionary equations with unbounded operator
coefficients in Banach spaces began to develop actively after the well-known monograph by
Henry \cite{Hen} (1981), the papers of Sacker, Sell \cite{Sacker_Sell1} (1994), the above mentioned mathematicians (Latushkin, Tomilov) and many others. Recently it has become clear about the importance of such research and in more general topological spaces such as the Frechet space (relevant results can be found in the work of Boichuk, Pokutnyi, Zhuravlyov \cite{BoiPokZhuravl} (2018) and the
earlier papers of Boichuk, Pokutnyi \cite{Boi_Pok}, \cite{Boichuk_Pokutnyi} (2014), where such problems were considered for the first time). In the paper Aragao \cite{Costa} (2019), the author gives an example pf a semigroup bounded operators $\{e^{m\Delta} : m \in \mathbb{N}\}$, where $\Delta$ is the Laplace operator in the unbounded domain, which is exponentially-dichotomous in the Frechet space, but not in the Banach space). In Aragao \cite{Costa} (2019), the author gives an example of a semigroup of bounded exponentially dichotomous in Frechet space, but not in Banach space. In the works of Walther \cite{Walt}, \cite{Walt1} (2016), the equation with delay in the space $C( ( - \infty, 0];\mathbb{R}^n)$ was investigated, which is a Frechet space (earlier, the concept of exponential dichotomy was used in the study of a system of differential equations with delay in the monograph \cite{BoiZhur} by Samoilenko, Boichuk, and Zhuravlev). The advantage in considering such a space is that it contains all the (complete) histories $x_{t} = x(t + \cdot), t \in \mathbb{R}$ of each solution of the differential equation $x'(t) = f(x_{t})$ in contrast from the Banach space. The exponential dichotomy for operators in the space of distributions also requires the use of Frechet spaces instead of Banach spaces.

Another method for studying the property of exponential dichotomy is the use of quadratic forms and sign-changing Lyapunov functions. In this formulation, this question for bounded matrix functions $A(t)$ was developed for the first time in the works of Malkin \cite{Malk} (1937) and then developed in the works of Daletskiy, Krein \cite{dal}, Valeev, Finin \cite{Val} (1980), Mitropolsky, Samoilenko, Kulik \cite{MitrSamKul} (1990), Stanzhitsky \cite{Stang} (2009) (the concept of exponential dichotomy was proposed in the study of Ito stochastic systems) and many others.

The concept of exponential dichotomy occupies a central place in the qualitative theory of
differential equations and dynamical systems. Systems with this property represent a class
whose trajectories can both increase and decrease at an exponential rate. An important property
is the roughness of the dichotomy for such systems, that is, the preservation of the dichotomy property under linear and nonlinear perturbations (which was discussed above, where the works of Massera and Schaeffer were mentioned, in which the corresponding roughness
property already appears implicitly [40]). Such problems are relevant and well developed for
finite-dimensional systems of differential equations, operator-differential equations with
bounded operator coefficients in the works of Daletsky, Krein \cite{dal} (1970) (using the principle
of contraction mappings and estimates for perturbing terms) and for equations with unbounded
operator coefficients in the monograph by Henry \cite{Hen} (1981) (approach using estimates for the evolution operator of the corresponding equation). Naulin, Pinto \cite{Naul} (1998), Vinograd \cite{Vinogr} (1988, 1991), Wiggins \cite{Ju} (2001), Popescu \cite{Pop}, \cite{Pop1} (2001, 2006), Zhou \cite{Lu} (2013), Lupa 
study of the roughness of various types of dichotomy (non-uniform, $(h, k)$-dichotomies, $\mu, \nu$ - dichotomies \cite{Lupa}, the approach proposed by Daletsky and Krein \cite{dal} is applied and further developed in the dichotomies, and so on) for operator-differential and difference equations in Banach spaces under linear perturbations. In the papers of Zang and Zong \cite{Zhang}, \cite{Zhang1} (2010, 2013), a similar method is used to study equations on time scales (including in the case of a parametric dependence of the perturbing term, the so-called robustness). In these papers, the proposed approach is extended to the case of unbounded operators, when the constants that define the dichotomy can be different at the corresponding intervals. Recently, the ideas of Zhu, Lu, Zang found its reflection in the study of the corresponding question for irreversible linear difference equations by Batelli, Frank and Palmer \cite{Bat} (2021). In \cite{Palmers} Palmer (1984) the roughness of the dichotomy is established on the semi-axe $\mathbb{R}_{+}$ provided that depending on time perturbation $B(t)$ tends to infinity. Chow \cite{Chow} (1995) studied the roughness for semiflows in Banach spaces. Feketa \cite{Feketa} (2015) studies the roughness of the dichotomy in the finite-dimensional case under the commutation condition of the summand $A(t)$ and perturbations using the Magnus expansion. Potzshe \cite{Potz} (2015) studied the preservation of the dichotomy property for difference equations under smoothly parametrically dependent linear perturbations. Note that in most papers, the roughness of the dichotomy under linear dichotomous on the entire axis $\mathbb{R}$ or semi-axes $\mathbb{R}_{+}$, $\mathbb{R}_{-}$. The main difference of this paper is that the roughness of the dichotomy is studied under the assumption that each of subsystems is exponentially-dichotomous. 

The main goal of this paper is to study the roughness of the dichotomy for a interconnected system of operator-differential equations with unbounded coefficients in the main part under both linear and nonlinear perturbations. The work has the following structure.

The first part contains the problem statement. The next two parts are devoted to obtaining sufficient conditions for the existence of generalized bounded solutions for a system of coupled equations and the corresponding estimates for the norms of solutions, as well as roughness conditions under linear perturbations by two methods, which are a development of the ideas of Daletsky, Krein \cite{dal}, Henry \cite{Hen} and Valeev \cite{Val},  which allow us to estimate solutions and prove roughnesss results. In Section 1 we introduce the basic concepts and notation. Section 2 contains the main results, with proofs, on exponential dichotomies for linear systems with unbounded operators. In Section 3 we indicate some extensions involving quadratic functionals. The concluding Section 4 discusses extensions to weakly nonlinear systems.
In contrast to earlier works, we directly study equations on the entire real line $\mathbb{R}$, not just on $\mathbb{R}_{+}$ or $\mathbb{R}_{-}$

{\bf 1. Statement of the problem. Linear case.} We consider interconnected system of differential equations with  stationar coefficients
\begin{equation} \label{ur:1}
x'(t) = (A + B)x(t), ~~t \in J
\end{equation}
where $x(t) = (x_1(t), x_2(t), ..., x_n(t))^{T} \in \mathbf{B} = B_{1}\times B_2 \times ... \times B_{n}$ is a vector valued function with values in the Banach space $\mathbf{B}$; $B_{i}, i = \overline{1, n}$ are Banach spaces; the linear operators $A$ and $B$ have the block matrix structure:
$$
A = \left( \begin{array}{ccccc} A_{11} & 0 & ... & 0 \\
0 & A_{22} & ... & 0 \\
... & ... & ... & ... \\
0 & 0 & ... & A_{nn}
\end{array}\right),~~
B = \left(\begin{array}{ccccc}
0 & A_{12} & A_{13} & ... & A_{1n} \\
A_{21} & 0 & A_{23} & ... & A_{2n} \\
... & ... & ... & ... & ... \\
A_{n1} & A_{n2} & ... & A_{nn} & 0
\end{array}\right).
$$
The operators $A_{ii}: B_{i} \rightarrow B_{i}$ are closed linear operators with domains $D(A_{ii}) \subset B_{i}$ dense in $B_{i}$, and the operators  $A_{ij} : B_{i} \rightarrow B_{j}, i \neq j$ are bounded: $A_{ij} \in \mathcal{L}(B_{i}, B_{j})$, $J = \mathbb{R}, \mathbb{R}_{+}$ or $\mathbb{R}_{-}$.
The coordinate form of system (\ref{ur:1}) has the following form:
\begin{equation} \label{eq:2}
\left\{ \begin{array}{ccccc} \frac{dx_1(t)}{dt} = A_{11}x_1(t) + A_{12}x_2(t) + ... + A_{1n}x_n(t), \\
\frac{dx_2(t)}{dt} = A_{21}x_1(t) + A_{22}x_2(t) + ... + A_{2n}x_n(t), \\
... \\
\frac{dx_n(t)}{dt} = A_{n1}x_1(t) + A_{n2}x_2(t) + ... + A_{nn}x_n(t).
\end{array} \right.
\end{equation}
The homogeneous (unperturbed) system has the form
\begin{equation} \label{ur:2}
x'(t) = Ax(t),
\end{equation}
or in coordinate form does not involve coupling between the components (with operators $A_{ii}$ which act in the corresponding spaces $B_{i}$)
\begin{equation} \label{eq:4}
\left\{ \begin{array}{ccccc}
\frac{dx_{1}(t)}{dt} = A_{11}x_{1}(t), \\
\frac{dx_{2}(t)}{dt} = A_{22}x_{2}(t), \\
... \\
\frac{dx_{n}(t)}{dt} = A_{nn}x_{n}(t).
\end{array}\right.
\end{equation}

Assume that each equation in system (\ref{eq:4}) (components of system (\ref{ur:2})) has an exponential dichotomy on $J$ (see \cite{dal}, \cite[p.245]{Hen}). It means that there exist projectors $P_{i}^{2} = P_{i}$, $M_{i}, N_{i}$, $\alpha_i, \beta_{i}$, $i = \overline{1,n}$ such that:
\begin{equation} \label{eq:5}
\left\{ \begin{array}{cccccc}
||e^{(t - s)A_{ii}}P_i|| \leq M_{i}e^{-\alpha_i(t - s)}, ~~t \geq s, \\
||e^{(t - s)A_{ii}}(I - P_{i})|| \leq N_{i}e^{-\beta_i(s - t)}, ~~t \leq s.
\end{array}\right.
\end{equation}
Here $e^{tA_{ii}}$ is the evolution operator of the $i$-th component (\ref{eq:4}) of the homogeneous system (\ref{ur:2}).
Moreover, the space $B_{i}$ can be decomposed in the direct sum of subspaces:
$$
B_{i} = P_{i}B_{i} \oplus (I - P_{i})B_{i}
$$
and the entire space
$$
\mathbf{B} = P_{-}\mathbf{B} \oplus P_{+}\mathbf{B} = \mathbf{B}_{-} \oplus \mathbf{B}_{+}
$$
with the projectors
$$
P_{-} = \left(\begin{array}{ccccc} P_{1} & 0 & ... & 0 \\
0 & P_{2} & ... & 0 \\
... & ... & ... & P_{n}
\end{array} \right), P_{+} = I - P_{-}
$$
and corresponding decomposition of the space $\mathbf{B}$ into invariant subspaces for the homogeneous system (\ref{ur:2}) (with bounded solutions of the unperturbed system (\ref{ur:2}) appearing from the corresponding subspaces).
Bounded solutions of the initial value problem (\ref{eq:5}) (see \cite{KrSavch}) 
$$
\frac{dx_{i}(t)}{dt} = A_{ii}x_{i}(t), ~~x_{i}(s) = x_{i 0} \in D(A_{ii})\cap P_{i}B_{i},
$$
satisfy the following estimate
$$
||x_{i}(t)|| \leq M_{i}e^{-\alpha_{i}(t - s)}||x_{i}(s)|| = M_{i}e^{-\alpha_{i}(t - s)}||x_{i0}||, ~~t \geq s
$$
for $x_{i0} \in P_{i}B_{i}$, with analogous estimate for $x_{i0} \in (I - P_{i})B_{i}$.
We find conditions on the perturbation terms $B$ (the coefficients $A_{ij}, i \neq j$) in system (\ref{ur:1}) which imply the existence of an exponential dichotomy for solutions of (\ref{ur:1}).

{\bf 2. Method of bounded solutions.} 1) a. At first we consider the case of half-lines $\mathbb{R}_{+}$, $\mathbb{R}_{-}$ for simplicity, then extend to the entire line $\mathbb{R}$.

If we will see on the perturbation term $Bx(t)$ as corresponding nonhomogeneity $f(t) = Bx(t)$, then we can show that the set of generalized bounded solutions on $\mathbb{R}_{+}$ of nonhomogeneous equaton (\ref{ur:1}) can be represented in the following form
\begin{equation} \label{ur:3}
x(t) = e^{tA}c_{-} + \int_{0}^{+\infty}G_{A}(t - \tau)f(\tau)d\tau,
\end{equation}
where \cite[p.119]{dal}, \cite{Rad}
$$
G_{A}(t) = \left\{ \begin{array}{ccccc} e^{tA}P_{-}, t > 0 \\
- e^{-tA}P_{+}, t < 0
\end{array} \right.
$$
is the Green's function for system (\ref{ur:2}). Here
$$
e^{tA} = \left( \begin{array}{ccccc}
e^{tA_{11}} & 0 & ... & 0\\
0 & e^{tA_{22}} & ... & 0 \\
... & ... & ... & e^{tA_{nn}}
\end{array} \right),
$$
$c_{-}$ is the initial value from the invariant subspace $\mathbf{B}_{-}$ so that $c_{-} = P_{-}c_{-}$.

Defining the operator $S_1$ on $BC(\mathbb{R}_{+}, \mathbf{B})$ (with the norm $|||f||| = \sup_{t \in \mathbb{R_{+}}}||f(t)||$) by the formula
\small{\begin{equation} \label{ur:5}
S_1x(t) = S_1\left( \begin{array}{ccccc}
x_1(t) \\
x_2(t) \\
... \\
x_{n}(t)
\end{array} \right):= \int_{0}^{+\infty}G_{A}(t - \tau)\left( \begin{array}{ccccc}
A_{12}x_2(\tau) + A_{13}x_3(\tau) + ... + A_{1n}x_n(\tau)\\
A_{21}x_1(\tau) + A_{23}x_3(\tau) + ... + A_{2n}x_n(\tau) \\
... \\
A_{n1}x_1(\tau) + A_{n2}x_2(\tau) +  ... + A_{nn-1}x_{n - 1}(\tau)
\end{array} \right)d\tau,
\end{equation}
}
we will show that under certain conditions this operator maps $BC(\mathbb{R}_{+}, \mathbf{B})$ into itself and is contractive (in such way we obtain the system of integral equations instead of differential system (\ref{eq:2})). Indeed, we have
$$
||S_1x(t)|| \leq \left(\int_{0}^{t}||G_{A}(t - \tau)||d\tau + \int_{t}^{+\infty}||G_{A}(t - \tau)||d\tau \right)\sum_{i, j = 1, i \neq j}^{n}||A_{ij}|||||x|||,
$$
or
$$
||S_1x(t)|| \leq \sqrt{n(n - 1)}\left(\int_{0}^{t}||G_{A}(t - \tau)||d\tau + \int_{t}^{+\infty}||G_{A}(t - \tau)||d\tau \right)\max_{i,j, i \neq j}\{||A_{ij}||\} |||x|||.
$$
Using the properties of the evolution operators we obtain the following inequalities:
$$
\int_{0}^{t}||G_{A}(t - \tau)||d\tau + \int_{t}^{+\infty}||G_{A}(t - \tau)||d\tau \leq \sum_{i = 1}^{n}\int_{0}^{t} ||e^{(t - \tau)A_{ii}}P_{i}|| d\tau +
$$
$$
+ \sum_{i = 1}^{n} \int_{t}^{+\infty} ||e^{-(t - \tau)A_{ii}}(I - P_{i})|| d\tau \leq \sum_{i = 1}^{n}M_{i}\int_{0}^{t} e^{-\alpha_i(t - \tau)} d\tau +
$$
$$
+ \sum_{i = 1}^{n}N_{i}\int_{t}^{+\infty} e^{-\beta_i(\tau - t)} d\tau = \sum_{i = 1}^{n}\left( \frac{M_i}{\alpha_i} + \frac{N_{i}}{\beta_{i}} - \frac{M_{i}}{\alpha_{i}}e^{-\alpha_{i}t}\right), ~~~t \geq 0
$$
or
$$
\int_{0}^{t}||G_{A}(t - \tau)||d\tau + \int_{t}^{+\infty}||G_{A}(t - \tau)||d\tau \leq \sqrt{n}\int_{0}^{t} \max_{i}\{||e^{(t - \tau)A_{ii}}P_{i}||\}d\tau +
$$
$$
+ \sqrt{n}\int_{t}^{+\infty} \max_{i}\{||e^{-(t - \tau)A_{ii}}(I - P_{i})||\}d\tau \leq \sqrt{n} \left(K_{1} + K_{2} - K_{1}e^{-\alpha t}\right), ~~~ t \geq 0
$$
where
$$
K_{1} = \max_{i = \overline{1,n}}\left\{ \frac{M_i}{\alpha_i}\right\}, K_{2} = \max_{i = \overline{1,n}}\left\{ \frac{N_i}{\beta_i}\right\}, \alpha = \min_{i = \overline{1, n}}\{ \alpha_i \}.
$$
From these estimates we obtain
\begin{equation} \label{nerivn:1}
|||S_1x||| = \sup_{t \in \mathbb{R}_{+}}||S_1x(t)|| \leq \sum_{i = 1}^{n}\left(\frac{M_i}{\alpha_i} + \frac{N_{i}}{\beta_{i}}\right) \sum_{i, j = 1, i \neq j}^{n}||A_{ij}|| |||x|||,
\end{equation}
or 
\begin{equation} \label{nerivn:2}
|||S_1x||| \leq n\sqrt{n - 1}\left(K_1 + K_2\right) \max_{i,j, i \neq j}\{||A_{ij}||\}|||x|||.
\end{equation}
Under the conditions
\begin{equation} \label{dost:1}
\sum_{i, j = 1, i \neq j}^{n}||A_{ij}|| < \frac{1}{\sum_{i = 1}^{n}\left(\frac{M_i}{\alpha_i} + \frac{N_{i}}{\beta_{i}}\right)}
\end{equation}
or
\begin{equation} \label{dost:2}
\max_{i,j, i \neq j}\{||A_{ij}||\} < \frac{1}{n\sqrt{n - 1}\left(K_1 + K_2\right)}
\end{equation}
the operator $S_1$ is contractive. Hence (\ref{eq:2}) has a unique solution given by the series
\begin{equation} \label{predst:1}
x(t) = (I - S_1)^{-1}e^{tA}c_{-}.
\end{equation}
This proves the following.

{\bf Lemma 1.} (Sufficient conditions for exponential dichotomy on $\mathbb{R}_{+}$). {\it Under conditions (\ref{dost:1}), (\ref{dost:2}) for any $c_{-} \in \mathbf{B}_{-}$ system (\ref{ur:1}) has a unique solution $x(t)$ which satisfies the following estimate
\begin{equation} \label{dostatochn:1}
|||x||| \leq \frac{\sum_{i = 1}^{n}M_{i}}{1 - \sum_{i = 1}^{n}\left(\frac{M_i}{\alpha_i} + \frac{N_{i}}{\beta_{i}}\right)\sum_{i,j = 1, j \neq j}^{n}||A_{ij}||}||c_{-}||
\end{equation}
or
\begin{equation} \label{dostatochn:2}
|||x||| \leq \frac{M}{1 - \sqrt{n(n - 1)}\left(K_1 + K_2\right)\max_{i,j, i \neq j}\{||A_{ij}||\}}||c_{-}||,
\end{equation}
where $M = \max_{i = \overline{1, n}}\{ M_i \}.$
}

{\bf Remark 1.} {\it From the obtained theorem as a corollary we have Theorem 4.2 in \cite[p.122]{dal}.}
In this case for the space of initial data we can obtain another decomposition:
$$
\mathbf{B} = \widetilde{P}_{-}\mathbf{B} \oplus \widetilde{P}_{+}\mathbf{B}
$$
where subspace $\widetilde{P}_{-}\mathbf{B}$ consists from the initial data in the following form: $c_{-} = x_{0} = \left(x_{10}, x_{20}, ..., x_{n0}\right)$ with corresponding bounded on the semi-axe $\mathbb{R}_{+}$ solutions. Using (\ref{predst:1}) we can obtain that for the space of initial data  the following representation is correct
$$
x_{0} = x(0) = Fc_{-}: = c_{-} + \int_{0}^{+ \infty}G_{A}(-\tau)Bx(\tau)d\tau = (I - P_{+}ZP_{-})c_{-},
$$
where
$$
Z = P_{+}ZP_{-} = \int_{0}^{+\infty}P_{+}e^{-\tau A}B(I - S_{1})^{-1}e^{\tau A}P_{-}d\tau
$$
is an operator with the following estimates 
$$
||Z|| \leq \frac{\sum_{i = 1}^{n}M_i \sum_{i = 1}^{n}\frac{N_{i}}{\beta_{i}}||B||}{1 - \sum_{i = 1}^{n}\left(\frac{M_i}{\alpha_i} + \frac{N_{i}}{\beta_{i}}\right)||B||}
\leq \frac{\sum_{i = 1}^{n}M_i\sum_{i = 1}^{n}\frac{N_{i}}{\beta_i}\sum_{i,j = 1, i \neq j}^{n}||A_{ij}||}{1 - \sum_{i = 1}^{n}\left(\frac{M_i}{\alpha_i} + \frac{N_{i}}{\beta_{i}}\right)\sum_{i,j = 1, i \neq j}^{n}||A_{ij}||}
$$
or
$$
||Z|| \leq \frac{\sqrt{n}M K_{2} \max_{i,j, i \neq j}\{||A_{ij}||\}}{1 - \sqrt{n(n - 1)}\left(K_1 + K_2\right)\max_{i,j, i \neq j}\{||A_{ij}||\}}.
$$
Hence the operator $F$ maps the subspace $P_{-}\mathbf{B}$ onto $\widetilde{P}_{-}\mathbf{B}$ with the estimate
\begin{equation} \label{fin:100}
\widetilde{P}_{-} = FP_{-}F^{-1}
\end{equation}
which carries out dichotomy of the perturbed system. Here $F^{-1} = I + P_{+}ZP_{-}$. 
Let's estimate parameters ($(\widetilde{M}, \mu$) of the perturbed system  (4). For any generalized bounded solution of the perturbed system (4) from the subspace 
$\widetilde{\mathbf{B}} \cap D(A) = \widetilde{P}_{-}\mathbf{B} \cap D(A)$ the following representation is hold
$$
x(t) = e^{(t - s)A}P_{-}x(s) + \int_{s}^{+\infty}G_{A}(t - \tau)Bx(\tau)d\tau,~ t \geq s.
$$
Then we can obtain the following estimates:
$$
||x(t)|| \leq \sum_{i = 1}^{n}M_{i}e^{-\alpha_{i}(t - s)}||x(s)|| + \sum_{i,j = 1, i \neq j}^{n}||A_{ij}||\sum_{i = 1}^{n}M_{i}\int_{s}^{t}e^{-\alpha_{i}(t - \tau)}||x(\tau)||d\tau  +
$$
\begin{equation}\label{nerav:1}
+ \sum_{i,j = 1, i \neq j}^{n}||A_{ij}||\sum_{i = 1}^{n}N_{i}\int_{t}^{+\infty}e^{-\beta_i(\tau - t)}||x(\tau)||d\tau
\end{equation}
or
$$
||x(t)|| \leq \sqrt{n}Me^{-\alpha(t - s)}||x(s)|| + \sqrt{n(n - 1)}\max_{i,j, i \neq j}\{||A_{ij}||\}\int_{s}^{t}\max_{i}\left\{M_{i}e^{-\alpha_{i}(t - \tau)}\right\}||x(\tau)||d\tau +
$$
\begin{equation} \label{nerav:2}
+ \sqrt{n(n - 1)}\max_{i,j, i \neq j}\{||A_{ij}||\}\int_{t}^{+\infty}\max_{i}\left\{N_{i}e^{-\beta_{i}(\tau - t)}\right\}||x(\tau)||d\tau.
\end{equation}
For the following kernels in (\ref{nerav:1}), (\ref{nerav:2}) 
$$
\mathcal{K}_1(t, \tau) = \left\{ \begin{array}{ccccc} \sum_{i,j = 1, i \neq j}^{n}||A_{ij}||\sum_{i = 1}^{n}M_{i}e^{-\alpha_{i}(t - \tau)}, \tau \leq t, \\
\sum_{i,j = 1, i \neq j}^{n}||A_{ij}||\sum_{i = 1}^{n}N_{i}e^{-\beta_{i}(\tau - t)}, \tau \geq t, \end{array} \right.
$$
and
$$
\mathcal{K}_2(t, \tau) = \left\{ \begin{array}{ccccc} \sqrt{n(n - 1)}\max_{i,j, i \neq j}\{||A_{ij}||\}\max_{i}\left\{M_ie^{-\alpha_{i}(t - \tau)}\right\}, \tau \leq t \\
\sqrt{n(n - 1)}\max_{i,j, i \neq j}\{||A_{ij}||\}\max_{i}\left\{ N_{i}e^{-\beta_{i}(\tau - t)}\right\}, \tau \geq t
\end{array} \right.
$$
respectively, the following estimates are hold
$$
\int_{s}^{+\infty}\mathcal{K}_1(t, \tau)d\tau \leq \sum_{i,j = 1, i \neq j}^{n}||A_{ij}||\sum_{i = 1}^{n}\left(\frac{M_{i}}{\alpha_{i}}\left( 1 - e^{-\alpha_{i}(t - s)}\right) + \frac{N_{i}}{\beta_{i}}\right) \leq
$$
$$
\leq \sum_{i,j = 1, i \neq j}^{n}||A_{ij}||\sum_{i = 1}^{n}\left(\frac{M_{i}}{\alpha_{i}} + \frac{N_{i}}{\beta_{i}}\right) < 1,
$$
$$
\int_{s}^{+\infty}\mathcal{K}_2(t, \tau)d\tau \leq \sqrt{n(n - 1)}\max_{i,j, i \neq j}\{||A_{ij}||\}\left(K_{1}\left( 1 - e^{-\alpha(t - s)}\right) + K_{2}\right) \leq
$$
$$
\leq \sqrt{n(n - 1)}\sum_{i,j = 1, i \neq j}^{n}||A_{ij}||\sum_{i = 1}^{n}\left(K_{1} + K_{2}\right) < 1
$$
by virtue of (13), (14) for all $t \in \mathbb{R}_{+}, t \geq s$. For simplicity of calculations, we strengthen the inequalities in (20), (21) as follows
$$
||x(t)|| \leq \sum_{i = 1}^{n}M_{i}e^{-\Lambda(t - s)}||x(s)|| +  \sum_{i,j = 1, i \neq j}^{n}||A_{ij}||\sum_{i = 1}^{n}M_{i}\int_{s}^{t}e^{-\Lambda(t - \tau)}||x(\tau)||d\tau +
$$
$$
+ \sum_{i,j = 1, i \neq j}^{n}||A_{ij}||\sum_{i = 1}^{n}N_{i}\int_{t}^{+\infty}e^{-\Lambda(\tau - t)}||x(\tau)||d\tau \leq \sum_{i = 1}^{n}M_{i}e^{-\Lambda(t - s)}||x(s)|| +
$$
\begin{equation} \label{nerav:3}
+  \overline{M}\sum_{i,j = 1, i \neq j}^{n}||A_{ij}||\int_{s}^{+\infty}e^{-\Lambda|t - \tau|}||x(\tau)||d\tau,
\end{equation}
$$
||x(t)|| \leq \sqrt{n}Me^{-\Lambda(t - s)}||x(s)|| + \sqrt{n(n - 1)}\max_{i,j, i \neq j}\{||A_{ij}||\}\max_{i}\left\{M_{i}\right\}\int_{s}^{t}e^{-\Lambda(t - \tau)}||x(\tau)||d\tau +
$$
$$
+ \sqrt{n(n - 1)}\max_{i,j, i \neq j}\{||A_{ij}||\}\max_{i}\left\{N_{i}\right\}\int_{t}^{+\infty}e^{-\Lambda(\tau - t)}||x(\tau)||d\tau \leq \sqrt{n}Me^{-\Lambda(t - s)}||x(s)|| +
$$
\begin{equation} \label{nerav:4}
+ \sqrt{n(n - 1)}\overline{N}\max_{i,j, i \neq j}\{||A_{ij}||\}\int_{s}^{+\infty}e^{-\Lambda|t - \tau|}||x(\tau)||d\tau,
\end{equation}
respectively. Here 
\begin{equation} \label{dost:3}
\Lambda = \min_{i = \overline{1, n}}\left\{ \alpha_{i}, \beta_{i} \right\},~ \overline{M} = \max \left\{ \sum_{i = 1}^{n}M_{i}, \sum_{i = 1}^{n}N_{i}\right\},~ \overline{N} = \max_{i = \overline{1, n}} \left\{ M_{i}, N_{i} \right\}.
\end{equation}
Applying corollary 2.3 of Lemma 2.1 \cite{dal}[p. 156] to the inequalities (22), (23) we can obtain the following estimates for $x(t)$ and assertion.

{\bf Theorem 1. } (sufficient condition of dichotomy).
{\it Every solution $x(t)$ of system (\ref{ur:1}) satisfies the following estimates:
$$
||x(t)|| \leq \widetilde{M}e^{-\mu(t - s)}||x(s)||, ~~~ t \geq s
$$
with constants
\begin{equation} \label{fin:1001}
\widetilde{M} = \frac{\sum_{i = 1}^{n}M_{i}}{\overline{M}\sum_{i,j = 1, i \neq j}^{n}||A_{ij}||}\Lambda,
\end{equation}
$$
\mu < \sqrt{\Lambda^{2} - 2\Lambda \overline{M} \sum_{i,j = 1, i \neq j}^{n}||A_{ij}||},
$$
under the condition
$$
\sum_{i,j = 1, i \neq j}^{n}||A_{ij}|| < \frac{\Lambda}{2\overline{M}}
$$
or
\begin{equation} \label{fin:1002}
\widetilde{M} = \frac{M\Lambda}{\sqrt{n - 1}\overline{N}\max_{i,j, i \neq j}\{||A_{ij}||\}},
\end{equation}
$$
\mu < \sqrt{\Lambda^{2} - 2\Lambda\sqrt{n(n - 1)}\overline{N}\max_{i,j, i \neq j}\{||A_{ij}||\}},
$$
under the condition
$$
\max_{i,j, i \neq j}\{||A_{ij}||\} < \frac{\Lambda}{2\sqrt{n(n - 1)}\overline{N}}.
$$
The constants $\Lambda, \overline{M}, \overline{N}$ are defined in (\ref{dost:3}).
}

{\bf Remark 2.} {\it To obtain another estimates on $\mu$ we can use the following technique. Let $\varphi(t) = ||x(t)||$. Setting $\varphi(t) = u(t)e^{-\mu(t - s)}$ ($\mu < \alpha_{i}, i = \overline{1, n}$) for $u(t)$ we have the following estimates
$$
u(t) \leq \sum_{i = 1}^{n}M_{i}e^{(-\alpha_i + \mu)(t - s)}u(s)  + \sum_{i,j = 1, i \neq j}^{n}||A_{ij}||\sum_{i = 1}^{n}M_{i}\int_{s}^{t}e^{-\alpha_{i}(t - \tau) + \mu (t - \tau)}u(\tau)d\tau  +
$$
$$
+ \sum_{i,j = 1, i \neq j}^{n}||A_{ij}||\sum_{i = 1}^{n}N_{i}\int_{t}^{+\infty}e^{-\beta_i(\tau - t) + \mu(t - \tau)}u(\tau)d\tau,
$$
or
$$
u(t) \leq  \sqrt{n}Me^{-(\alpha - \mu)(t - s)}u(s) + \sqrt{n(n - 1)}\max_{i,j, i \neq j}\{||A_{ij}||\}M\int_{s}^{t}e^{-\alpha(t - \tau) + \mu(t - \tau)}u(\tau)d\tau +
$$
$$
+ \sqrt{n(n - 1)}\max_{i,j, i \neq j}\{||A_{ij}||\}N\int_{t}^{+\infty}e^{-\beta(\tau - t) + \mu(t - \tau)}u(\tau)d\tau
$$
respectively. Here $\beta = \min_{i = \overline{1, n}}\{ \beta_{i} \} $, $N = \max_{i = \overline{1, n}}\{ N_{i} \}$.
Consider the following operator $J$
$$
(Ju)(t): = \sum_{i,j = 1, i \neq j}^{n}||A_{ij}||\sum_{i = 1}^{n}M_{i}\int_{s}^{t}e^{-\alpha_{i}(t - \tau) + \mu (t - \tau)}u(\tau)d\tau
$$
$$
+ \sum_{i,j = 1, i \neq j}^{n}||A_{ij}||\sum_{i = 1}^{n}N_{i}\int_{t}^{+\infty}e^{-\beta_i(\tau - t) + \mu(t - \tau)}u(\tau)d\tau
$$
or
$$
(Ju)(t) = \sqrt{n(n - 1)}\max_{i,j, i \neq j}\{||A_{ij}||\}M\int_{s}^{t}e^{-\alpha(t - \tau) + \mu(t - \tau)}u(\tau)d\tau +
$$
$$
+ \sqrt{n(n - 1)}\max_{i,j, i \neq j}\{||A_{ij}||\}N\int_{t}^{+\infty}e^{-\beta(\tau - t) + \mu(t - \tau)}u(\tau)d\tau
$$
respectively, with the sup norm on $[s; +\infty)$ given by $||u|| = \sup_{t \in [s; +\infty)}|u(t)|$ (we consider the space of continuous functions).
It is easy to see that
$$
||J|| \leq \sum_{i,j = 1, i \neq j}^{n}||A_{ij}||\sum_{i = 1}^{n}\left(\frac{M_{i}}{\alpha_{i} - \mu} + \frac{N_{i}}{\beta_{i} + \mu} - \frac{M_{i}}{\alpha_{i} - \mu}e^{(\mu - \alpha_{i})(t - s)}\right) \leq
$$
$$
\leq \overline{N} \max_{i}\{ \alpha_{i} + \beta_{i} \} \sum_{i,j = 1, i \neq j}^{n}||A_{ij}|| \sum_{i = 1}^{n} \frac{1}{(\alpha_i - \mu)(\beta_i + \mu)} \leq \frac{n \overline{N} \max_{i}\{ \alpha_{i} + \beta_{i} \}\sum_{i,j = 1, i \neq j}^{n}||A_{ij}|| }{\Lambda^{2} - \mu^{2}},
$$
or
$$
||J|| \leq \sqrt{n(n - 1)}\max_{i,j, i\neq j}\{||A_{ij}||\}\left( \frac{M - Me^{(- \alpha + \mu)(t - s)}}{\alpha - \mu} + \frac{N}{\beta + \mu}\right) \leq
$$
$$
\leq \frac{\sqrt{n(n - 1)}\overline{N}\max_{i,j, i\neq j}\{||A_{ij}||\}}{\Lambda^{2} - \mu^{2}}
$$
respectively. Condition $||J|| < 1$ guarantees that the operator $J$ is contractive and we can obtain estimates of the parameter $\mu$.
In such way we obtain the following estimates
$$
||x(t)|| \leq \frac{\sum_{i = 1}^{n}M_{i}}{1 - ||J||}e^{-\mu(t - s)}||x(s)||, ~~ t \geq s
$$
or
$$
||x(t)|| \leq \frac{\sqrt{n}M}{1 - ||J||}e^{-\mu(t - s)}||x(s)||, ~~ t \geq s.
$$


Setting
$$
\frac{n \overline{N} \max_{i}\{ \alpha_{i} + \beta_{i} \}\sum_{i,j = 1, i \neq j}^{n}||A_{ij}|| }{\Lambda^{2} - \mu^{2}} = q_1,
$$
and
$$
\frac{\sqrt{n(n - 1)}\overline{N}\max_{i,j, i\neq j}\{||A_{ij}||\}}{\Lambda^{2} - \mu^{2}} = q_2
$$
respectively, we obtain the following theorem.
}

{\bf Theorem 2.} (sufficient condition of dichotomy). {\it Every solution $x(t)$ of system (\ref{ur:1}) satisfies the following estimate:
$$
||x(t)|| \leq \widetilde{M}e^{-\mu(t - s)}||x(s)||, ~~ t \geq s,
$$
with constants $\widetilde{M}$ given by the formulas
\begin{equation} \label{fin:1000}
\widetilde{M} = \frac{\sum_{i = 1}^{n}M_{i}}{1 - q_1}, ~~\mbox{or}~~ \widetilde{M} = \frac{\sqrt{n}M}{1 - q_2},
\end{equation}
and parameter $\mu$ satisfy the inequalities
{\small \begin{equation} \label{dost:4}
\mu < \sqrt{\Lambda^{2} - n\overline{N}\max_{i = \overline{1, n}}\{\alpha_{i} +\beta_{i}\}\sum_{i,j = 1, j\neq j}^{n}||A_{ij}||}, ~~\mbox{under}~~ \sum_{i, j = 1, i \neq j}^{n}||A_{ij}|| < \frac{\Lambda^{2}}{n\overline{N}\max_{i = \overline{1, n}}\{ \alpha_{i} + \beta_{i}\}}
\end{equation}
}
or
\begin{equation} \label{dost:5}
\mu < \sqrt{\Lambda^{2} - \sqrt{n(n - 1)}\overline{N}\max_{i,j, i \neq j}\{||A_{ij}||\}}, ~~\mbox{under}~~ \max_{i, j, i\neq j}\{||A_{ij}||\} < \frac{\Lambda^{2}}{\sqrt{n(n - 1)}\overline{N}}
\end{equation}
respectively.
}

On the half-line $\mathbb{R}_{-}$ we consider the following system of integral equations
\begin{equation} \label{ura:1}
\left( \begin{array}{ccccc}
x_1(t) \\
x_2(t) \\
... \\
x_{n}(t)
\end{array} \right)= e^{tA}c_{+} + \int_{-\infty}^{0}G_{A}(t - \tau)\left( \begin{array}{ccccc}
A_{12}x_2(\tau) + ... + A_{1n}x_n(\tau)\\
A_{21}x_1(\tau) + ... + A_{2n}x_n(\tau) \\
... \\
A_{n1}x_1(\tau) + ... + A_{nn-1}x_{n - 1}(\tau)
\end{array} \right)d\tau,
\end{equation}
where $c_{+} = P_{+}c_{+} \in \mathbf{B}_{+}$, and the operator
\begin{equation} \label{ur:5}
S_2x(t) = S_2\left( \begin{array}{ccccc}
x_1(t) \\
x_2(t) \\
... \\
x_{n}(t)
\end{array} \right):= \int_{-\infty}^{0}G_{A}(t - \tau)\left( \begin{array}{ccccc}
A_{12}x_2(\tau) + ... + A_{1n}x_n(\tau)\\
A_{21}x_1(\tau) + ... + A_{2n}x_n(\tau) \\
... \\
A_{n1}x_1(\tau) + ... + A_{nn-1}x_{n-1}(\tau)
\end{array} \right)d\tau.
\end{equation}

Similar estimates as for $\mathbb{R}_{+}$ show that the operator $S_2$ maps $BC(\mathbb{R}_{-}, \mathbf{B})$ into itself with estimates analogous to (\ref{nerivn:1}), (\ref{nerivn:2}) for $S_1$:
\begin{equation} \label{urf:1}
|||S_2x||| = \sup_{t \in \mathbb{R}_{-}}||S_2x(t)|| \leq \sum_{i = 1}^{n}\left(\frac{M_i}{\alpha_i} + \frac{N_{i}}{\beta_{i}}\right)\sum_{i,j = 1, j \neq j}^{n}||A_{ij}|||||x|||,
\end{equation}
or
\begin{equation} \label{urf:2}
|||S_2x||| \leq 2\sqrt{n(n - 1)}\left(K_1 + K_2\right) \max_{i,j, i \neq j}\{||A_{ij}||\}|||x|||
\end{equation}
respectively. In such way we obtain the following assertion.

{\bf Lemma 2.} (Sufficient conditions of exponential dichotomy on $\mathbb{R}_{-}$).
{\it Under conditions (\ref{urf:1}), (\ref{urf:2}) for any $c_{+} \in \mathbf{B}_{+}$ system (\ref{ur:1}) has a unique generalized solution $x(t)$ which satisfies the following estimates 
\begin{equation} \label{urf:3}
|||x||| \leq \frac{\sum_{i = 1}^{n}N_{i}}{1 - \sum_{i = 1}^{n}\left( \frac{M_i}{\alpha_i} + \frac{N_{i}}{\beta_{i}}\right)\sum_{i,j = 1, j \neq j}^{n}||A_{ij}||}||c_{+}||
\end{equation}
or
\begin{equation} \label{urf:4}
|||x||| \leq \frac{N}{1 - 2\sqrt{n(n - 1)}\left(K_1 + K_2\right)\max_{i,j, i \neq j}\{||A_{ij}||\}}||c_{+}||
\end{equation}
respectively.
Here $N = \max_{i = \overline{1, n}}\{ N_{i} \}$.
}

In this case for the construction of a bounded solution consider initial conditions of the form $c_{+} = x_{0}$ with bounded initial values on $\mathbb{R}_{-}$. Formula (\ref{ura:1}) shows that for such initial data the solution is representable as
$$
x(0) = Hc_{+}: = c_{+} + \int_{-\infty}^{0}G_{A}(-\tau)Bx(\tau)d\tau = (I - P_{-}Z'P_{+})c_{+},
$$
where
$$
Z' = P_{-}Z'P_{+} = \int_{-\infty}^{0}P_{-}e^{-\tau A}B(I - S_{2})^{-1}e^{\tau A}P_{+}d\tau.
$$
For the operator $Z'$ estimates similar to the $\mathbb{R}_{+}$ case hold:
$$
||Z'|| \leq \frac{\sum_{i = 1}^{n}N_i \sum_{i = 1}^{n}\frac{M_{i}}{\alpha_{i}}||B||}{1 - \sum_{i = 1}^{n}\left(\frac{M_i}{\alpha_i} + \frac{N_{i}}{\beta_{i}}\right)||B||}
\leq \frac{\sum_{i = 1}^{n}N_i\sum_{i = 1}^{n}\frac{M_{i}}{\alpha_i}\sum_{i,j = 1, i \neq j}^{n}||A_{ij}||}{1 - \sum_{i = 1}^{n}\left(\frac{M_i}{\alpha_i} + \frac{N_{i}}{\beta_{i}}\right)\sum_{i,j = 1, i \neq j}^{n}||A_{ij}||}
$$
or
$$
||Z'|| \leq \frac{\sqrt{n}N K_{1} \max_{i,j, i \neq j}\{||A_{ij}||\}}{1 - \sqrt{n(n - 1)}\left(K_1 + K_2\right)\max_{i,j, i \neq j}\{||A_{ij}||\}}.
$$

Hence the operator $H$ maps the subspace $P_{+}\mathbf{B}$ onto $\widetilde{P}_{+}\mathbf{B}$ with projectors
$$
\widetilde{P}_{+} = HP_{+}H^{-1}
$$
which carries out dochotomy of the perturbed system (4). Here $H^{-1} =I + P_{-}Z'P_{+}$.

1) b. For the entire line $\mathbb{R}$ and an each nonhomogeneity $f(t) \in BC(\mathbb{R}, \mathbf{B})$  there exists unique bounded on $\mathbb{R}$ generalized solution  (\ref{eq:1}),
(corresponding unperturbed system (\ref{ur:2} has only trivial bounded solution):
\begin{equation} \label{neodn:2}
x(t) = S x(t) = S\left( \begin{array}{ccccc}
x_1(t) \\
x_2(t) \\
... \\
x_{n}(t)
\end{array} \right):= \int_{-\infty}^{+\infty}G_{A}(t - \tau)\left( \begin{array}{ccccc}
A_{12}x_2(\tau) + ... + A_{1n}x_n(\tau)\\
A_{21}x_1(\tau) + ... + A_{2n}x_n(\tau) \\
... \\
A_{n1}x_1(\tau) + ... + A_{nn-1}x_{n-1}(\tau)
\end{array} \right)d\tau.
\end{equation}
Estimates for the operator $S$ have the following form:
$$
||Sx(t)|| \leq \sum_{i = 1}^{n}\left(\frac{M_i}{\alpha_i} + \frac{N_{i}}{\beta_{i}}\right)\sum_{i, j = 1, j \neq j}^{n}||A_{ij}|||||x|||, ~~~t \in \mathbb{R}
$$
or
$$
||Sx(t)|| \leq  \sqrt{n(n - 1)}\left(K_1 + K_2\right) \max_{i,j, i \neq j}\{||A_{ij}||\}|||x|||, ~~~t \in \mathbb{R},
$$
and hence in the uniform metric have the same form as for the operators using properties of $S_1$ and $S_2$ ((\ref{nerivn:1}), (\ref{nerivn:2}), (\ref{urf:1}), (\ref{urf:2})).
These properties imply a uniqueness result for solutions of (\ref{nel:1}).

In this case the required invariant splitting has the form
$$
\mathbf{B} = \widetilde{P}_{-}\mathbf{B}\oplus \widetilde{P}_{+}\mathbf{B} = LP_{-}\mathbf{B} \oplus LP_{+}\mathbf{B}
$$
with $L = FP_{-} + HP_{+}$, or equivalently,
$$
L = I - P_{+}ZP_{-} - P_{-}Z'P_{+} = I - Z - Z'
$$
as can be verified. For uniqueness it suffices that
$$
||Z + Z'|| < 1,
$$
which holds under the conditions
\begin{equation} \label{osn:1}
\sum_{i,j = 1, j \neq j}^{n}||A_{ij}|| \leq \frac{1}{\sum_{i = 1}^{n}N_{i}\sum_{i = 1}^{n}\frac{M_{i}}{\alpha_{i}} + \sum_{i = 1}^{n}M_{i}\sum_{i = 1}^{n}\frac{N_{i}}{\beta_{i}} + \sum_{i = 1}^{n}\left(\frac{M_{i}}{\alpha_{i}} + \frac{N_{i}}{\beta_{i}}\right)}
\end{equation}
or
\begin{equation} \label{osn:2}
\max_{i,j, i \neq j}\{||A_{ij}||\} \leq \frac{1}{\sqrt{n}(NK_1 + MK_2) + \sqrt{n(n - 1)}(K_1 + K_2)}.
\end{equation}

This proves the existence and uniqueness of solutions for system (\ref{ur:1}) under conditions of Theorems 1, 2 and Corollaries 1, 2.

{\bf 3. Method of Lyapunov functions.} We find sufficient conditions of exponential dichotomy for the system (\ref{ur:1}), (\ref{eq:2})  in the case when $\mathbf{B} = \mathbf{H} = \mathcal{H}_1 \times \mathcal{H}_2 \times ... \times \mathcal{H}_n$, $\mathcal{H}_{i}, i = \overline{1, n}$ are Hilbert spaces with using of quadratic forms.

We will use different criterion with using of Lyapunov functions and quadratic forms (see \cite{dal}, \cite[p.58, Theorem 2.3]{Val}, \cite{MitrSamKul}).

By virtue of dichotomy, there exist such operators $C_{i}, i = \overline{1, n}$  satisfying the following Lyapunov equations
\begin{equation} \label{eq:7}
A_{ii}^{*}C_i + C_{i}A_{ii} = - H_i,~~ i = \overline{1, n}
\end{equation}
with positive definite operators $H_i, i = \overline{1, n}$. We can consider quasidiagonal form of operators $H_{i0}$ 
$$
H_{i0} = P_{i}^{*}H_iP_{i} + (I - P_{i})^{*}H_i(I - P_{i}), ~~ i = \overline{1, n}.
$$
for simplicity.
In the case when operators $H_{i}, i = \overline{1, n}$ are identity operators, we get the set of operators
$$P_{i}^{*}P_{i} + (I - P_{i})^{*}(I - P_{i}),~~ i = \overline{1, n},$$
respectively. 

In the case of dichotomy on the whole axis for each subsystem (\ref{eq:7}) it is sufficient to find operators $C_{i}, i = \overline{1, n}$ in the following form (see \cite[p.60]{Val}, \cite[p.60, Theorem 7.1]{dal} by analogy):
\begin{equation} \label{eq:8}
C_{i} = \int_{0}^{+\infty}e^{tA_{ii}^{*}}P_{i}^{*}P_{i}e^{tA_{ii}}dt - \int_{-\infty}^{0}e^{tA_{ii}^{*}}(I - P_{i})^{*}(I - P_{i})e^{tA_{ii}}dt.
\end{equation}
The following estimates are hold
\begin{equation} \label{eq:10}
||C_{i}|| \leq \int_{0}^{+\infty}||e^{tA_{ii}^{*}}P_{i}^{*}||||P_ie^{tA_{ii}}||dt + \int_{-\infty}^{0}||e^{tA_{ii}^{*}}(I - P_{i})^{*}||||(I - P_{i})e^{tA_{ii}}|| dt.
\end{equation}
By virtue of (\ref{eq:5}) we can obtain the following inequalities
\begin{equation} \label{eq:11}
||C_{i}|| \leq M_{i}^{2}\int_{0}^{+\infty}e^{-2\alpha_{i} t}dt + N_{i}^{2}\int_{- \infty}^{0} e^{2\alpha_it}dt = \frac{M_{i}^{2}}{2\alpha_i} + \frac{N_{i}^{2}}{2\beta_{i}}.
\end{equation}
Therefore for the operator $A$ of system (\ref{ur:2}) operator Lyapunov matrix
$$
C = \left( \begin{array}{ccccc} C_{1} & 0 & ... & 0\\
0 & C_{2} & ... & 0 \\
... & ... & ... & ... \\
0 & ... & 0 & C_{n}
\end{array}\right)
$$
satisfies the operator Lyapunov equation
\begin{equation} \label{eqs:1}
A^{*}C + CA = - H,
\end{equation}
where
\begin{equation} \label{eqs:2}
H = \left(\begin{array}{cccccc}
 P_{1}^{*}P_1 + (I - P_{1})^{*}(I - P_{1}) & 0 & ... & 0 \\
 ... & ... & ... & ... \\
 0 & ... & ... & P_{n}^{*}P_n + (I - P_{n})^{*}(I - P_{n})
\end{array} \right)
\end{equation}
is a positive definite operator.

Consider the quadratic functional $v(X) = (CX, X)$ with the operator matrix $C$ satisfying (\ref{eqs:1}) with negative-definite derivative by virtue of (\ref{ur:2}).

If $P_{i} = P_{i}^{*}, i = \overline{1, n}$, then $H = I$. Computing the derivative we obtain
$$
\frac{dv(X)}{dt} = \left((A^{*} + B^{*})CX, X\right) + \left((C(A + B)X, X\right) =
$$
$$
= -(HX, X) + (B^{*}CX, X) + (CBX, X) = -(X,X) + (B^{*}CX, X) + (CBX, X).
$$
Under the condition
\begin{equation} \label{wer:1}
||C||(||B^{*}|| + ||B||) < 1
\end{equation}
the functional $v(X)$ is a positive definite quadratic form (it follows from Valeev). From (\ref{eq:11}) and Minkovsky inequaility follows that 
$$
||C|| \leq \sum_{i = 1}^{n}||C_{i}|| \leq \sum_{i = 1}^{n}\left(\frac{M_{i}^{2}}{2\alpha_i} + \frac{N_{i}^{2}}{2\beta_i}\right)
$$
or
$$
||C|| \leq \sqrt{n}\max_{i = \overline{1, n}}||C_{i}|| \leq \sqrt{n} \max_{i = \overline{1, n}} \left\{ \frac{M_{i}^{2}}{2\alpha_i} + \frac{N_{i}^{2}}{2\beta_i}\right\}.
$$
Direct estimates give
$$
||B|| + ||B^{*}|| \leq \sum_{i,j = 1, i \neq j}^{n}\left(||A_{ij}|| + ||A_{ij}^{*}||\right)   =  2\sum_{i,j = 1, i \neq j}^{n}||A_{ij}||
$$
or
$$
||B|| + ||B^{*}|| \leq \sqrt{n(n - 1)}\left(\max_{i,j, i \neq j}\{||A_{ij}||\} + \max_{i,j, i \neq j}\{||A_{ij}^{*}||\}\right) = 2\sqrt{n(n - 1)}\max_{i,j, i \neq j}\{||A_{ij}||\}
$$
respectively. Substituting into (\ref{wer:1}) proves positivity under the smallness conditions.

{\bf Theorem 3.} {\it If $P_{i} = P_{i}^{*}, i = \overline{1, n}$, then under the conditions:
\begin{equation} \label{osnov:1}
\sum_{i,j = 1, i \neq j}^{n}||A_{ij}|| < \frac{1}{\sum_{i = 1}^{n}\left(\frac{M_{i}^{2}}{\alpha_i} + \frac{N_{i}^{2}}{\beta_{i}}\right)},
\end{equation}
or
\begin{equation} \label{osnov:2}
\max_{i, j, i \neq j}\{||A_{ij}||\} < \frac{1}{n\sqrt{n - 1}\max_{i = \overline{1, n}} \left\{ \frac{M_{i}^{2}}{\alpha_{i}} + \frac{N_{i}^{2}}{\beta_i}\right\}}
\end{equation}
the functional $(CX, X)$ is positive definite, and perturbed system (\ref{ur:1}) has exponential dichotomy.
}

{\bf Remark 3.} {\it In general case it is sufficient that $B^{*}C + CB << H$ (i.e. $A >> 0$ is a positive definite operator).}

{\bf 4. Nonlinear case.} Finally consider nonlinear perturbed interconnected system of the form
\begin{equation} \label{nel:1}
x'(t) = Ax(t) + R(x(t)),~~ t \in \mathbb{R}
\end{equation}
where block structure of unbounded oeprator $A$ has the same form as in the linear case, and the nonlinear term $R$ has the structure:
$$
R(x(t)) = \left( \begin{array}{ccccc} R_{1}(x_{2}(t), x_{3}(t), ..., x_{n}(t)) \\
R_{2}(x_{1}(t), x_{3}(t), ..., x_{n}(t)) \\
... \\
R_{n}(x_{1}(t), x_{2}(t), ..., x_{n - 1}(t))
\end{array} \right),
$$
where nonlinearities
$$
R_{i}(\cdot, ..., \cdot): B_{1}\times B_{2} \times ... \times B_{i - 1} \times B_{i + 1} \times ... \times B_{n} \rightarrow B_{i}
$$
are locally Lipschitz continuous (with Lipschitz constants on bounded sets).

{\bf Definition 1.} {\it Write that
$\{ R \in (T_{i}, L_{i}, \rho), i = \overline{1, n}\}$ on set $S_{\rho}$ (here $S_{\rho} = \{x = (x_{1}, x_{2}, ..., x_{n})|:$ $~~ ||x||_{\mathbf{B}} \leq \rho \}$) if:

1. $||R_{i}(x_{1}, ..., x_{i - 1}, x_{i + 1}, ..., x_{n})|| \leq T_{i}$;

2. For any $x^{1} = (x_{1}^{1}, ..., x_{i}^{1}, ..., x_{n}^{1}), ~x^{2} = (x_{1}^{2}, ..., x_{i}^{2}, ..., x_{n}^{2}) \in S_{\rho}$
$$
||R_{i}(x_{1}^{2}, ..., x_{i - 1}^{2}, x_{i + 1}^{2}, ..., x_{n}^{2}) - R_{i}(x_{1}^{1}, ..., x_{i - 1}^{1}, x_{i + 1}^{1}, ..., x_{n}^{1})|| \leq L_{i}\sqrt{\sum_{j = 1, j \neq i}^{n}||x_{j}^{2} - x_{j}^{1}||^{2}};
$$

3. $R_{i}(0, ..., 0) = 0$, $i = \overline{1, n}$.
}

By analogy with Section 2 we will prove the existence of bounded generalized solution of the nonlinear equation (\ref{nel:1}) in the following form
\begin{equation} \label{neodn:2}
y(t) = S x(t) = S\left( \begin{array}{ccccc}
x_1(t) \\
x_2(t) \\
... \\
x_{n}(t)
\end{array} \right):= \int_{-\infty}^{+\infty}G_{A}(t - \tau)\left( \begin{array}{ccccc}
R_{1}(x_{2}(\tau), x_{3}(\tau), ..., x_{n}(\tau))\\
R_{2}(x_{1}(\tau), x_{3}(\tau), ..., x_{n}(\tau)) \\
... \\
R_{n}(x_{1}(\tau), x_{2}(\tau), ..., x_{n - 1}(\tau))
\end{array} \right)d\tau
\end{equation}
which for each fixed $t \in \mathbb{R}$ belongs to the set $S_{\rho}$.
Really, for $x(t) \in S_{\rho}$ with bounded norm $|||x||| = \sup_{t \in \mathbb{R}}||x(t)|| \leq \rho$ we obtain that:
$$
|||y||| \leq \sum_{i = 1}^{n}\left(\frac{M_i}{\alpha_i} + \frac{N_{i}}{\beta_{i}}\right) \sum_{i = 1}^{n}|||R_{i}||| \leq \sum_{i = 1}^{n}\left(\frac{M_i}{\alpha_i} + \frac{N_{i}}{\beta_{i}}\right) \sum_{i = 1}^{n}T_{i}  \leq \rho
$$
or
$$
|||y||| \leq \sqrt{n(n - 1)}\left(K_1 + K_2\right) \max_{i = \overline{1, n}} T_{i} \leq \rho
$$
and needed condition we can obtain by choosing parameters $T_{i}$:
\begin{equation} \label{nonlinear:1}
\sum_{i = 1}^{n}T_{i} \leq \frac{\rho}{\sum_{i = 1}^{n}\left(\frac{M_i}{\alpha_i} + \frac{N_{i}}{\beta_{i}}\right)}
\end{equation}
or
\begin{equation} \label{nonlinear:2}
\max_{i} T_{i} \leq \frac{\rho}{\sqrt{n(n - 1)}\left(K_1 + K_2\right)}.
\end{equation}
Next we obtain the following estimates:
\begin{equation} \label{nonlinear:3}
|||Sx^{2} - Sx^{1}||| \leq \sum_{i = 1}^{n}\left(\frac{M_i}{\alpha_i} + \frac{N_{i}}{\beta_{i}}\right) \sum_{i = 1}^{n}L_{i} |||x^{2} - x^{1}|||
\end{equation}
or
\begin{equation} \label{nonlinear:4}
|||Sx^{2} - Sx^{1}||| \leq \sqrt{n(n - 1)}\left(K_1 + K_2\right) \max_{i = \overline{1, n}} L_{i} |||x^{2} - x^{1}|||.
\end{equation}
From the estimates (\ref{nonlinear:1})-(\ref{nonlinear:4}) follows the following assertion.

{\bf Lemma 3.} {\it For any $\rho > 0$ there are constants $T_{i}, L_{i}, i = \overline{1, n}$ satisfying (\ref{nonlinear:1}), (\ref{nonlinear:2}) with  $R \in (T_{i}, L_{i}, \rho)$ such that under the following conditions
$$
\sum_{i = 1}^{n}L_{i} < \frac{1}{\sum_{i = 1}^{n}\left(\frac{M_i}{\alpha_i} + \frac{N_{i}}{\beta_{i}}\right)}
$$
or
$$
\max_{i = \overline{1, n}} L_{i} < \frac{1}{\sqrt{n(n - 1)}\left(K_1 + K_2\right)}
$$
then system (\ref{nel:1}) has a unique solution $x(t)$ with
$$
|||x||| = \sup_{t \in \mathbb{R}}||x(t)|| \leq \rho.
$$
}

By analogy with linear case we can write generalized bounded solution of the nonlinear system (\ref{nel:1})  for $t \geq s$  in the following form
$$
x(t) = e^{A(t - s)}x(s) + \int_{s}^{+\infty}G_A(t - \tau)R(x(\tau))d\tau,
$$
with $x(s) = P_{-}x(s)$.

We find conditions such that trajectory $x(t)$ belongs to $S_{\rho}$:
$$
||x(t)|| \leq \sum_{i = 1}^{n}M_{i}e^{-\alpha_{i}(t - s)}||x(s)|| + \sum_{i = 1}^{n}\left(\frac{M_i}{\alpha_i}  + \frac{N_{i}}{\beta_i} - \frac{M_{i}}{\alpha_i}e^{-\alpha_i(t - s)}\right)\sum_{i = 1}^{n}T_{i} \leq
$$
$$
\leq \sum_{i = 1}^{n}M_{i}||x(s)|| + \sum_{i = 1}^{n}\left(\frac{M_i}{\alpha_i}  + \frac{N_{i}}{\beta_i}\right)\sum_{i = 1}^{n}T_{i}
$$
or
$$
||x(t)|| \leq \sqrt{n}Me^{-\alpha(t - s)}||x(s)|| + \sqrt{n(n - 1)}\left(K_1 + K_2 - K_1 e^{-\alpha(t - s)} \right)\max_{i = \overline{1, n}}T_{i} \leq
$$
$$
\leq \sqrt{n}M||x(s)|| + \sqrt{n(n - 1)}\left(K_1 + K_2\right)\max_{i = \overline{1, n}}T_{i}.
$$
If
$$
x(s) \in \mathbf{B}_{-} \cap S_{\frac{\rho}{2\sum_{i = 1}^{n}M_{i}}} ~\mbox{or}~ x(s) \in \mathbf{B}_{-} \cap S_{\frac{\rho}{2\sqrt{n}M}},
$$
then under the following conditions on parameters $T_{i}$
\begin{equation} \label{nonlinear:5}
\sum_{i = 1}^{n}T_{i} \leq \frac{\rho}{2\sum_{i = 1}^{n}\left(\frac{M_i}{\alpha_i}  + \frac{N_{i}}{\beta_i}\right)}
\end{equation}
or
\begin{equation} \label{nonlinear:6}
\max_{i = \overline{1, n}}T_{i} \leq \frac{\rho}{2\sqrt{n(n - 1)}\left(K_1 + K_2\right)}
\end{equation}
obtain what we needed.

For estimating of parameters of dichotomy for $t \geq s$ we can use the following inequalities:
$$
||x(t)|| \leq \sum_{i = 1}^{n}M_{i}e^{-\alpha_{i}(t - s)}||x(s)||  +  \sum_{i = 1}^{n}L_{i}\sum_{i = 1}^{n}M_{i}\int_{s}^{t}e^{-\alpha_{i}(t - \tau)}||x(\tau)||d\tau  +
$$
\begin{equation} \label{nonlinear:8}
+ \sum_{i = 1}^{n}L_{i}\sum_{i = 1}^{n}N_{i}\int_{t}^{+\infty}e^{-\beta_i(\tau - t)}||x(\tau)||d\tau
\end{equation}
or
$$
||x(t)|| \leq \sqrt{n}Me^{-\alpha(t - s)}||x(s)|| + \sqrt{n(n - 1)}\max_{i}L_{i}\int_{s}^{t}\max_{i}\left\{M_{i}e^{-\alpha_{i}(t - \tau)}\right\}||x(\tau)||d\tau +
$$
\begin{equation} \label{nonlinear:9}
+ \sqrt{n(n - 1)}\max_{i}L_{i}\int_{t}^{+\infty}\max_{i}\left\{N_{i}e^{-\beta_{i}(\tau - t)}\right\}||x(\tau)||d\tau.
\end{equation}
Applying Theorems 1, 2 we obtain the following assertions.

{\bf Theorem 4.} (sufficient condition of dichotomy). {\it For any $\rho > 0$ and constants $L_{i}, T_{i}$ satisfying (\ref{nonlinear:5}), (\ref{nonlinear:6}) with $R \in (T_{i}, L_{i}, \rho)$, if $x(s) \in \mathbf{B}_{-} \cap S_{\frac{\rho}{2\sum_{i = 1}^{n}M_{i}}}$ or $x(s) \in \mathbf{B}_{-} \cap S_{\frac{\rho}{2\sqrt{n}M}}$ then the solution satisfies
the following inequality
$$
||x(t)|| \leq \widetilde{M}e^{-\mu(t - s)}||x(s)||, ~~t \geq s
$$
with parametr 
$$
\widetilde{M} = \frac{\sum_{i = 1}^{n}M_{i}}{\overline{M} \sum_{i = 1}^{n}L_{i}}\Lambda, ~~ \mu < \sqrt{\Lambda^{2} - 2\Lambda \overline{M} \sum_{i = 1}^{n}L_{i}},
$$
under the condition
$$
 \sum_{i = 1}^{n}L_{i} < \frac{\Lambda}{2\overline{M}}
$$
or
$$
\widetilde{M} = \frac{M\Lambda}{\sqrt{n - 1}\overline{N}\max_{i}L_{i}},~~
\mu < \sqrt{\Lambda^{2} - 2\Lambda\sqrt{n(n - 1)}\overline{N}\max_{i}L_{i}},
$$
under the condition
$$
\max_{i}L_{i} < \frac{\Lambda}{2\sqrt{n(n - 1)}\overline{N}}.
$$
The constants $\Lambda, \overline{M}, \overline{N}$ are defined in (\ref{dost:3}). Moreover, the trajectory stays in the ball $S_{\rho}$.
}

{\bf Theorem 5.} (sufficient condition of dichotomy). {\it For any $\rho > 0$ and constants $L_{i}, T_{i}$ satisfying (\ref{nonlinear:5}), (\ref{nonlinear:6}) with $R \in (T_{i}, L_{i}, \rho)$, if $x(s) \in \mathbf{B}_{-} \cap S_{\frac{\rho}{2\sum_{i = 1}^{n}M_{i}}}$ or $x(s) \in \mathbf{B}_{-} \cap S_{\frac{\rho}{2\sqrt{n}M}}$ the solution satisfies the following estimate
$$
||x(t)|| \leq \widetilde{M}e^{-\mu(t - s)}||x(s)||,~~ t \geq s
$$
with parameters $\widetilde{M}$ given by
$$
\widetilde{M} = \frac{\sum_{i = 1}^{n}M_{i}}{1 - q_1}, ~~\mbox{or}~~ \widetilde{M} = \frac{\sqrt{n}M}{1 - q_2},
$$
and $\mu$ satisfying the inequalities
{\small $$
\mu < \sqrt{\Lambda^{2} - n\overline{N}\max_{i = \overline{1, n}}\{\alpha_{i} +\beta_{i}\}\sum_{i = 1}^{n}L_{i}}, ~~\mbox{under}~~ \sum_{i = 1}^{n}L_{i} < \frac{\Lambda^{2}}{n\overline{N}\max_{i = \overline{1, n}}\{ \alpha_{i} + \beta_{i}\}}
$$}
or
$$
\mu < \sqrt{\Lambda^{2} - \sqrt{n(n - 1)}\overline{N}\max_{i}L_{i}}, ~~\mbox{under}~~ \max_{i}L_{i} < \frac{\Lambda^{2}}{\sqrt{n(n - 1)}\overline{N}}
$$
respectively. Moreover, the trajectory stays in the ball $S_{\rho}$.
}

Consider a special case when the nonlinearity $R$ is Frechet differentiable and the following representation is hold
$$
R(x) = Bx + \widetilde{R}(x)
$$
or in coordinate form
\begin{equation} \label{fin:1}
R_{i}(x_{1}, ..., x_{i - 1}, x_{i + 1}, .., x_{n}) = \sum_{j = 1, j \neq j}^{n}A_{ij}x_{j} + \widetilde{R}_{i}(x_{1}, ..., x_{i - 1}, x_{i + 1}, ..., x_{n}),
\end{equation}
where the linear and bounded operators $A_{ij}: B_{j} \rightarrow A_{i}, A_{ij} \in \mathcal{L}(B_{j}, B_{i})$ are the same as under linear perturbations (\ref{ur:1}) and nonlinearities $\{ \widetilde{R} \in (T_{i}, L_{i}, \rho), i = \overline{1, n}\}$.
We consider the case of semi-axe $\mathbb{R}_{+}$ for simplicity. Using theorems 1, 2 we get conditions which guarantee that linear perturbed system (\ref{ur:1}) is exponentially dichotomous. Then there is projector $\widetilde{P}_{-}^{2} = \widetilde{P}_{-}$ (\ref{fin:100}), constants $\widetilde{M}_{1,2}, \mu$ such that for the evolution operator $e^{t(A + B)}$ of the perturbed system (\ref{ur:1}) the following estimates are hold:
\begin{equation} \label{fin:2}
\left\{ \begin{array}{ccccc} ||e^{(t - s)(A + B)}\widetilde{P}_{-}|| \leq \widetilde{M}_1e^{-\mu(t - s)}, ~~t \geq s \\
||e^{(t - s)(A + B)}(I - \widetilde{P}_{-})|| \leq \widetilde{M}_2e^{-\mu(s - t)}, ~~ t \leq s,
\end{array} \right.
\end{equation}
where $\widetilde{M}_{1}$ is given by (\ref{fin:1001}), (\ref{fin:1002}), (\ref{fin:1000}),  and $\widetilde{M}_2$ by the formulas
$$
\widetilde{M}_2 = \frac{\sum_{i = 1}^{n}N_{i}}{\overline{M}\sum_{i,j = 1, i \neq j}^{n}||A_{ij}||}\Lambda,~~
\widetilde{M}_2 = \frac{N\Lambda}{\sqrt{n - 1}\overline{N}\max_{i,j, i \neq j}\{||A_{ij}||\}},
$$
or
$$
\widetilde{M}_2 = \frac{\sum_{i = 1}^{n}N_{i}}{1 - q_1}, ~~ \widetilde{M}_2 = \frac{\sqrt{n}N}{1 - q_2}.
$$
Then we can prove existence of generalized bounded solution of the nonlinear system as a solution of the following system of integral equations
\begin{equation} \label{neodn:2}
y(t):= e^{t(A + B)}\widetilde{c}_{-} + \int_{0}^{+\infty}G_{A + B}(t - \tau)\left( \begin{array}{ccccc}
\widetilde{R}_{1}(x_{2}(\tau), x_{3}(\tau), ..., x_{n}(\tau))\\
\widetilde{R}_{2}(x_{1}(\tau), x_{3}(\tau), ..., x_{n}(\tau)) \\
... \\
\widetilde{R}_{n}(x_{1}(\tau), x_{2}(\tau), ..., x_{n - 1}(\tau))
\end{array} \right)d\tau,
\end{equation}
where $\widetilde{c}_{-} = \widetilde{P}_{-}\widetilde{c}_{-}$ which belongs to the ball  $S_{\rho}$ for any $t \in \mathbb{R}_{+}$ under the following conditions
$$
\sum_{i = 1}^{n}T_{i} \leq \frac{\rho \mu}{\widetilde{M}_1 + \widetilde{M}_2}, ~~\mbox{or}~~\max_{i}T_{i} \leq \frac{\rho \mu }{\widetilde{M}_1 + \widetilde{M}_2}.
$$
From the theorem 4, 5 we obtain the following assertion.

{\bf Corollary 4.} (sufficient condition of dichotomy). {\it For any $\rho > 0$ there are constants $L_{i}, T_{i}$ such that under the following conditions

1.$ \sum_{i = 1}^{n}T_{i} < \frac{\rho \mu}{\widetilde{M}_1 + \widetilde{M}_2},~$ or $~\max_{i}T_{i} < \frac{\rho \mu}{2\widetilde{M}_3}; $

2. $ \sum_{i = 1}^{n}L_i < \frac{\mu}{2\widetilde{M}_3}~,$ or $~\max_{i}L_{i} < \frac{\mu}{2\widetilde{M}_3};$

3. $\widetilde{R} \in (T_{i}, L_{i}, \rho)$, $\widetilde{M}_3 = \max\{ \widetilde{M}_1, \widetilde{M}_2 \}$,

for any $x(s) \in \widetilde{P}_{-}\mathbf{B} \cap S_{\frac{\rho}{2\left(\widetilde{M}_1 + \widetilde{M}_2\right)}}$  or $x(s) \in \widetilde{P}_{-}\mathbf{B} \cap S_{\frac{\rho}{4\widetilde{M}_3}}$ the following estimate is hold
$$
||x(t)|| \leq \widetilde{M}e^{-\nu(t - s)}||x(s)||, ~~t \geq s
$$
where parameter $\widetilde{M}$ given by
$$
\widetilde{M} = \frac{\widetilde{M}_1}{\widetilde{M}_3 \sum_{i = 1}^{n}L_i}\mu,~~ \widetilde{M} = \frac{\widetilde{M}_1}{\widetilde{M}_3 \max_{i}L_i}\mu,
$$
or
$$
\widetilde{M} = \frac{\widetilde{M}_1(\mu^{2} - \nu^2)}{\mu^{2} - \nu^{2} - 2\mu\sum_{i = 1}^{n}L_{i}\widetilde{M}_3},~~~ \widetilde{M} = \frac{\widetilde{M}_1(\mu^{2} - \nu^2)}{\mu^{2} - \nu^{2} - 2\mu\max_{i}L_{i}\widetilde{M}_3},
$$
and $\nu$ satisfies the following conditions
$$
\nu < \sqrt{\mu^2 - 2\mu\widetilde{M}_3 \sum_{i = 1}^{n}L_i}~~ \mbox{or}~~\nu < \sqrt{\mu^2 - 2\mu\widetilde{M}_3 \max_{i}L_i}.
$$
}

\section{Conclusions}

Proposed in the paper method gives possibility to apply it for the resonance problems (when the uniqueness solutions can disturbed. Moreover, in the case when the considered problem can be ill-posed). It should be noted then we can use this method for obtaining the conditions of the so-called input-to-state stability. Such theory originates from the papers of Sontag (see for example \cite{Sontag-1989} - \cite{Sontag-1998}) and then developed in the papers of the following authors (see \cite{Dashkovskiy-Mironchenko-2012-01} - \cite{Dashkovskiy-Slynko-Kapustyan-2021}, \cite{Jacob-Mironchenko-2019}, \cite{Mironchenko-Ito-2014} -  \cite{Mironchenko-Wirth-2017}).

\end{document}